\documentclass[a4paper,11pt]{article}

\usepackage{amssymb}

\newtheorem{lem}{Lemma}[section]
\newtheorem{prp}[lem]{Proposition}
\newtheorem{df}[lem]{Definition}
\newtheorem{thm}[lem]{Theorem}

\newlength{\eqdemoffset}
\newlength{\enumrqueoffset}
\newenvironment{dem}[1][]
               {\begin{trivlist}\item[]{\sc Proof#1.}~}
               {\ifhmode{\unskip~\nobreak}\else%
                {\nopagebreak\vspace{\eqdemoffset}\leavevmode\hfill}\fi%
                \hfill$\blacksquare$\end{trivlist}}
\newenvironment{rque}[1]
               {\begin{trivlist}\item[]{\sc #1.}~}
               {\ifhmode{\unskip~\nobreak}\else%
                {\nopagebreak\vspace{\enumrqueoffset}\leavevmode\hfill}\fi%
                \hfill$\square$\end{trivlist}}

\newcommand{\text}[1]{\mbox{#1}}
\newcommand{\verylongarrow}[1]{\stackrel{\displaystyle #1}%
  {\displaystyle{-\hspace{-0.5em}-\hspace{-0.5em}\longrightarrow}}}
\newcommand{\free}{*}
\newcommand{\red}{\mathrm{red}}
\newcommand{\Cst}{$C^*$-\relax}
\newcommand{\tens}{\makebox[.8em][c]{$\otimes$}}
\newcommand{\id}{{\rm id}}
\newcommand{\rond}{\makebox[1em][c]{$\circ$}}
\newcommand{\CC}{\mathbb C}
\newcommand{\ZZ}{\mathbb Z}
\newcommand{\NN}{\mathbb N}

\title{$K$-amenability for amalgamated free products of amenable discrete quantum groups}
\author{Roland {\sc Vergnioux}} 
\date{}

\begin{document}

\maketitle

\begin{abstract}
  The basic notions and results of equivariant $KK$-theory concerning crossed products can
  be extended to the case of locally compact quantum groups. We recall these constructions
  and prove some usefull properties of subgroups and amalgamated free products of discrete
  quantum groups. Using these properties and a quantum analogue of the Bass-Serre tree, we
  establish the $K$-amenability of amalgamated free products of amenable discrete quantum
  groups. \\[\smallskipamount]
\end{abstract}

The main goal of this paper is to prove the $K$-amenability of an amalgamated free product of
amenable discrete quantum groups. The notion of $K$-amenability was first introduced by Cuntz
\cite{Cuntz:kmoy} for discrete groups: the aim was to give a simpler proof to a result of
Pimsner and Voiculescu \cite{PimsnerVoicu:freecross} calculating the $K$-theory of the reduced
\Cst algebra of a free group. Cuntz proves that the $K$-theory of the reduced and full \Cst
algebras of the free groups are the same, and gives in \cite{Cuntz:freeprod} a very simple way
to compute it in the full case.

In \cite{JulgValette:kmoy}, Julg and Valette extend the notion of $K$-amenability to the
locally compact case and establish the $K$-amenability of locally compact groups acting on
trees with amenable stabilizers. By a construction of Serre \cite{Serres:arbres}, this includes
the case of amalgamated free products of amenable discrete groups. To prove the $K$-amenability
of a locally compact group $G$, one has to construct an element $\alpha \in KK_G(\CC,\CC)$
using representations of $G$ that are weakly contained in the regular one, and then to prove
that $\alpha$ is homotopic to the unit element of $KK_G(\CC,\CC)$. In \cite{JulgValette:kmoy},
both of the steps are carried out in a very geometrical way. Moreover, it turns out that
$\alpha$ can be interpreted as the $\gamma$ element used to prove the Baum-Connes conjecture in
this context \cite{KaspSkand:buildings_novikov}.

In fact, the $K$-amenability of the starting groups is often sufficient to get the
$K$-amenability of their free product. This was already the case in the original result of
Cuntz \cite{Cuntz:kmoy}, which deals with the non-amalgamated and non-quantum case.  This has
been generalized to the amalgamated and non-quantum case by Pimsner \cite{Pimsner:kktrees} (in
the wider framework of group actions on trees). Finally, the non-amalgamated and quantum case
is covered by the work of Germain \cite{Germain:freereduced}. In this paper we treat the
amalgamated and quantum case, but assume that both starting groups are amenable.

\bigskip

In the first section of this paper, we state some results that extend to the case of locally
compact quantum groups classical tools of equivariant $KK$-theory related with crossed
products, including the notion of $K$-amenability. The general framework is the equivariant
$KK$-theory with respect to coactions of Hopf \Cst algebras \cite{BaajSkand:kks} on the one
hand, and the theory of locally compact quantum groups \cite{BaajSkand:unit,KustVaes:lcqg} on
the other hand.

In the second section, we prove some useful results on subgroups and amalgamated free products
of discrete quantum groups. Here we base ourselves on the definitions given by Wang in
\cite{Wang:freeprod}, and extend some of the results therein. The main goal is to give an
alternative form for the regular representation of an amalgamated free product of discrete
quantum groups, which we then use for the $K$-theoretic construction of the last section.

Finally we construct in the third section the Julg-Valette operator associated to the ``quantum
tree'' of an amalgamated free product of amenable discrete quantum groups. We prove that it
defines an element of $KK$-theory which is homotopic to the unit element, thus establishing the
$K$-amenability of the amalgamated free product under consideration.

\section{Equivariant $KK$-theory for locally compact \\ quantum groups}

In fact, most of the work of extending Kasparov's equivariant theory to locally compact quantum
groups was done by Baaj and Skandalis in \cite{BaajSkand:kks}, where they define equivariant
$KK$-groups with respect to coactions of Hopf \Cst algebras, and prove the existence of the
Kasparov product in this setting. When the Hopf \Cst algebra under consideration is associated
to a locally compact quantum group, some additional features of the classical theory are
expected to extend to the quantum case, especially in relation with duality and crossed
products. It turns out that it is indeed what happens: we will give here a review of the
results obtained in \cite{Vergnioux:these}. We refer the interested reader to the same
reference for the proofs, which contain almost no difficulty particular to the quantum case.
Note that some of these results were already proved in \cite{BaajSkand:kks} in the case of
duals of locally compact groups.

\bigskip

The framework of our study will be the one of \Cst algebraic locally compact quantum groups, as
defined by Kustermans and Vaes in \cite{KustVaes:lcqg}. However, we mainly work with the Kac
systems $(H,V,U)$ associated to these quantum groups. Let us recall some notations and
definitions from \cite{BaajSkand:unit}.

Let $H$ be a Hilbert space, we will denote by $\Sigma$ the flip operator on $H\tens H$, and we
will use the ``leg numbering'' notation for tensor products --- for instance $\Sigma_{23} =
\id\tens\Sigma \in L(H\tens H\tens H)$. A unitary operator $V \in L(H\tens H)$ is said to be
multiplicative if it verifies the pentagon equation $V_{12}V_{13}V_{23} = V_{23}V_{12}$. We
define two closed subalgebras of $L(H)$ associated to $V$ in the following way:
\begin{eqnarray*}
  && S_\red = \overline{\mathrm{Lin}}~ \{ (\omega\tens\id)(V) ~|~ \omega\in L(H)_* \} \text{,} \\
  && \hat S_\red = \overline{\mathrm{Lin}}~ \{ (\id\tens\omega)(V) ~|~ \omega\in L(H)_* \} \text{.}
\end{eqnarray*}
More generally, if $X$ is a representation of $V$ on some Hilbert $B$-module $E$, i.e. a unitary
morphism $X\in L_B(E\tens H)$ such that $X_{12}X_{13}V_{23} = V_{23}X_{12}$, we put $\hat S_X =
\overline{\mathrm{Lin}}~ \{ (\id\tens\omega)(X) ~|~ \omega\in L(H)_* \} \subset L_B(E)$.
Besides, if $U$ is an involutive unitary on the same Hilbert space $H$, we put $\hat V = \Sigma
(U\tens 1) V (U\tens 1) \Sigma$ and $\tilde V = \Sigma(1\tens U)V(1\tens U)\Sigma$.  A Kac
system in the (strong) sense of \cite{BaajSkand:unit} is a triple $(H,V,U)$ such that
$(V(U\tens 1)\Sigma)^3 = 1$, $V$ and $\hat V$ are multiplicative, and both fulfill an
analytical condition called regularity.

The Kac systems $(H,V,U)$ associated to locally compact quantum groups are not Kac systems in
the sense of \cite{BaajSkand:unit}, however they verify the following useful properties
\cite{KustVaes:lcqg,BaajSkand:unit,Woro:multiplunit}:
\begin{enumerate}  \setlength{\itemsep}{0ex}
  \renewcommand\theenumi{{\it\roman{enumi}}} \renewcommand\labelenumi{{\it\roman{enumi}}.}
\item \label{item:kac_faible_mult} $V$ and $\hat V$ are multiplicative,
\item $[S_\red, US_\red U] = [\hat S_\red, U\hat S_\red U] = 0$,
\item \label{item:kac_faible_cstar} $S_\red$ and $\hat S_\red$ are \Cst algebras,
\item \label{item:kac_faible_repr} if $X$ is a representation of $V$ or $\hat V$, it belongs to
  $M(\hat S_X\tens K(H))$.
\end{enumerate}
In fact these are the only properties required for the validity of our results on equivariant
$KK$-theory, which means that we do not make a direct use of the Haar weights on the locally
compact quantum groups under consideration. Note that Properties~\ref{item:kac_faible_cstar}
and~\ref{item:kac_faible_repr} are automatically verified when $V$ and $\hat V$ are regular
\cite{BaajSkand:unit}, semi-regular \cite{Baaj:E(2)} or manageable \cite{Woro:multiplunit}.

\bigskip

If $(H,V,U)$ is the Kac system associated to a locally compact quantum group, and more
generally if it verifies properties \ref{item:kac_faible_mult}--\ref{item:kac_faible_repr}, the
\Cst algebras $S_\red$ and $\hat S_\red$ are naturally endowed with coproducts:
\begin{eqnarray*}
  &&\delta_\red(s) = V(s\tens 1)V^* = \hat V^* (1\tens s)\hat V \text{,} \\
  &&\hat\delta_\red(\hat s) = V^*(1\tens \hat s)V = \tilde V(\hat s\tens 1)\tilde V^* \text{.}
\end{eqnarray*}
Moreover one can check that $(S_\red,\delta_\red)$ and $(\hat S_\red,\hat\delta_\red)$ are then
Hopf \Cst algebras, in particular $\delta_\red$ is a non-degenerate homomorphism from $S_\red$
to $M(S_\red\tens S_\red)$ and verifies the coassociativity identity
$(\id\tens\delta_\red)\rond \delta_\red = (\delta_\red\tens\id)\rond \delta_\red$.

Recall that a coaction of $S_\red$ on a \Cst algebra $A$ is a non-degenerate homomorphism
$\delta_A : A \to M(A\tens S_\red)$ such that $\delta_A(A) (1\tens S_\red) \subset (A\tens
S_\red)$ and $(\id\tens\delta_\red) \rond \delta_A = (\delta_A\tens\id) \rond \delta_A$. We say
that $A$ is a $S_\red$-algebra if $\delta_A$ is injective and $\delta_A(A)(1\tens S_\red)$
spans a dense subspace of $A\tens S_\red$. We say that $a\in A$ is invariant if $\delta_A(a) =
a\tens 1$, and we denote by $A^{S_\red}$ the sub-\Cst algebra of invariant elements. There is
also a notion of coaction on Hilbert \Cst modules: a coaction of $S_\red$ on a Hilbert
$A$-module $E$ is a linear map $\delta_E : E \to M(E\tens S_\red)$ compatible with
$\delta_\red$, $\delta_A$ and the Hilbertian structure --- see \cite{BaajSkand:kks} for the
details.  When $\delta_A$ is trivial, i.e. $A^{S_\red} = A$, the coactions $\delta_E$ of $S_\red$
on $E$ correspond in fact to the representations $X$ of $V$ on $E$, via the following formula:
$\delta_E(\xi) = X(\xi\tens 1)$.

Finally, let $\delta_A$, $\delta_B$, $\delta_E$ be coactions of $S_\red$ on two \Cst algebras
$A$, $B$ and a Hilbert $B$-module $E$. We say that a representation $\pi : A \to L_B(E)$ is
covariant if
\begin{displaymath}
  \forall~ a\in A,~ \xi\in E~~
  \delta_E(\pi(a)\xi) = (\pi\tens\id)\delta_A(a) \cdot \delta_E(\xi) \text{.}
\end{displaymath}
If the coaction on $B$ is trivial and $X$ is the representation of $V$ associated to
$\delta_E$, the covariance of $\pi$ can also be written $(\pi\tens\id)\delta_A(a) = X
(\pi(a)\tens 1) X^*$.

\bigskip

Let $A$ be a \Cst algebra endowed with a coaction of $S_\red$. Let $X$ be a representation of
$V$ on a Hilbert $B$-module $E$, where $S_\red$ coacts trivially on $B$. To each covariant
representation $\pi : A \to L_B(E)$ one can associate a crossed product \Cst algebra in the
following way:
\begin{displaymath}
  A\rtimes_\pi\hat S = \overline{\mathrm{Lin}}~ \{ \pi(a) x ~|~ 
  a\in A,~~ x\in\hat S_X \} \text{.}
\end{displaymath}
In particular, there is a universal covariant representation $(\pi^A, X^A)$ giving a full
crossed product $A\rtimes\hat S$, and a regular covariant representation $(\pi^A_\red,
X^A_\red)$ giving a reduced crossed product $A\rtimes_\red\hat S$. Let us call $A_1$ the \Cst
algebra $A$ endowed with the trivial coaction of $S_\red$. The regular covariant representation
of $A$ is then defined to be $\delta_A : A \to M(A\tens K(H))$, with $E = A\tens H$ seen as a
$A_1$-module and equipped with the representation $X^A_\red = 1\tens V$.

One important tool is the dual coaction of $\hat S_\red$ (or of its full version $\hat S$) on
$A\rtimes_\red\hat S$ (or $A\rtimes\hat S$). For instance we use the following notations for
the coactions of $\hat S_\red$:
\begin{eqnarray*}
  &&\hat\delta_{A\rtimes\hat S} : A\rtimes\hat S \to M((A\rtimes\hat S)\tens\hat S_\red)
  \text{~~~and} \\ &&\hat\delta_{A\rtimes_\red\hat S} : A\rtimes_\red\hat S \to 
  M((A\rtimes_\red\hat S)\tens\hat S_\red) \text{.}
\end{eqnarray*}
These coactions are naturally related by the reduction homomorphism $\lambda_A : A\rtimes\hat S
\to A\rtimes_\red \hat S$, and they arise from the following weak inclusions of covariant
representations:
\begin{displaymath}
  (\pi^A\tens\id, X_{13}^A V_{23}) \prec (\pi^A, X^A) \text{~~~and~~~}
  (\pi^A_\red\tens\id, X_{\red,13}^A V_{23}) \prec (\pi^A_\red, X^A_\red) \text{.}
\end{displaymath}
There is another weak inclusion $(\pi^A\tens\id, X_{13}^A V_{23}) \prec (\pi^A_\red, X^A_\red)$
which induces a ``false coaction'' $\delta'_A : A\rtimes_\red\hat S \to (A\rtimes\hat S) \tens
\hat S_\red$. The homomorphism $\delta'_A$ is one of the important tools for the proof of
Theorem~\ref{thm:K_amenability} on $K$-amenability.

To conclude on crossed products, let us present a simple characterization of $A\rtimes_\red\hat
S$ when $S_\red$ is unital (i.e. the quantum group is compact) and $A$ is a $S_\red$-algebra. We
introduce a new coaction of $S_\red$ on $E = A\tens H$, seen as an $A$-module: we put
$\delta'_E(a\tens\xi) = ((U\tens 1)V(U\tens 1)(\xi\tens 1))_{23} \delta_A(a)_{13}$. This
coaction induces a coaction of $S_\red$ on $K_A(E)$, and we have then $A\rtimes_\red\hat S =
K_A(E)^{S_\red}$. This equality is the main tool for the proof of the ``Green-Julg''
Theorem~\ref{thm:green_julg}.

\bigskip

Now we turn to equivariant $KK$-theory. To simplify the statements we will assume that all the
\Cst algebras we consider are separable. Let $(H,V,U)$ be the Kac system associated to a
locally compact quantum group, and let $S_\red$ be its reduced Hopf \Cst algebra. Let
$\delta_B$ be a coaction of $S_\red$ on some \Cst algebra $B$. If $\delta_E$ is a coaction of
$S_\red$ on a Hilbert $B$-module $E$, we provide the $B$-module $E\tens H$ with the following
coaction of $S_\red$: $\delta_{E\tens H}(\zeta\tens\xi) = (V(\xi\tens 1))_{23}
\delta_E(\zeta)_{13}$. Let us call $\mathcal{H}_B$ the Hilbertian direct sum of countably many
copies of the Hilbert $B$-module $B$.  We have then the following equivariant ``Kasparov''
stabilization theorem:

\begin{thm}
  Let $(H,V,U)$ be the Kac system associated to a locally compact quantum group, and let
  $S_\red$ be its reduced Hopf \Cst algebra. Let $\delta_B$ be a coaction of $S_\red$ on a \Cst
  algebra $B$.
  \begin{enumerate}
  \item Let $E$, $F$ be two Hilbert $B$-modules equipped with coactions of $S_\red$. If $E$ and
    $F$ are isomorphic as Hilbert $B$-modules, then there is an equivariant isomorphism between
    $E\tens H$ and $F\tens H$.
  \item Let $E$ be a countably generated Hilbert $B$-module endowed with a coaction of
    $S_\red$.  Then $(E\oplus \mathcal{H}_B)\tens H$ is equivariantly isomorphic to
    $\mathcal{H}_B\tens H$.
  \item Assume moreover that $S_\red$ is unital. Then $E \oplus (\mathcal{H}_B\tens H)$ is
    equivariantly isomorphic to $\mathcal{H}_B\tens H$.
  \end{enumerate}
\end{thm}

Let us now recall the definition of $KK_{S_\red}(A,B)$, where the \Cst algebras $A$ and $B$ are
endowed with coactions of $S_\red$: we denote by $E_{S_\red}(A,B)$ the set of triples
$(E,\pi,F)$ where $E$ is a countably generated, $\ZZ/2\ZZ$-graded Hilbert $B$-module endowed
with a coaction of $S_\red$, $\pi$ is a covariant, grading preserving representation of $A$ on
$E$, and $F\in L_B(E)$ is a morphism of degree $1$ such that
\begin{eqnarray*} 
  &&[\pi(A), F] \text{,~~~} \pi(A) (F^2-1) \text{~~~and~~~}
  \pi(A) (F-F^*) ~\subset~ K_B(E) \text{,} \\
  &&\pi(A)\tens S_\red ~ (F\tens 1 - \delta_{K(E)} (F)) \subset K_{B\tens S_\red}(E\tens
  S_\red) \text{.}
\end{eqnarray*}
Then $KK_{S_\red}(A,B)$ is the quotient of $E_{S_\red}(A,B)$ with respect to the homotopy
relation induced by $E_{S_\red}(A, B[0,1])$. Like in the classical case, every
equivariant homomorphism $\phi : A\to B$ defines an element $[\phi] \in KK_{S_\red}(A,B)$.

Using the fact that the Hopf \Cst algebras under consideration are associated to some locally
compact quantum groups, one can define descent morphisms for these equivariant $KK$-groups.  If
$E$ is a Hilbert $B$-module endowed with a coaction of $S_\red$, the full and reduced crossed
products of $E$ by $S_\red$ are the relative tensor products $E\rtimes\hat S =
E\tens_B(B\rtimes\hat S)$ and $E\rtimes_\red\hat S = E\tens_B(B\rtimes_\red\hat S)$. One can
prove that the \Cst algebra $K(E\rtimes\hat S)$ naturally identifies with the crossed product
$K(E)\rtimes\hat S$, so that any \Cst representation $\pi : A \to M(K(E))$ gives rise to a
homomorphism $\pi\rtimes\id : A\rtimes\hat S \to M(K(E\rtimes\hat S))$. The same works for
reduced crossed products, and this indeed defines descent morphisms:
\begin{thm} --- cf \cite[rque 7.7b]{BaajSkand:unit}. 
  Let $S_\red$ be the reduced Hopf \Cst algebra associated to a locally compact quantum group,
  and $A$, $B$ two \Cst algebras equipped with coactions of $S_\red$. 
  \begin{enumerate}
  \item For any $(E,\pi,F) \in E_{S_\red}(A,B)$, the triple $(E\rtimes\hat S, \pi\rtimes\id,
    F\tens_B\id)$ is an element of $E(A\rtimes\hat S, B\rtimes\hat S)$, and this defines a
    group morphism $j : KK_{S_\red}(A,B) \rightarrow KK(A\rtimes\hat S, B\rtimes\hat S)$. 
  \item This morphism is compatible with the Kasparov product: for all $x\in KK_{S_\red}(A,D)$,
    $y\in KK_{S_\red}(D,B)$ one has $j(x\tens y) = j(x)\tens j(y)$.
  \item In the same way, considering reduced crossed products one obtains a group morphism
    $j_\red : KK_{S_\red}(A,B) \rightarrow KK(A\rtimes_\red\hat S, B\rtimes_\red\hat S)$
    compatible with the Kasparov product.
  \end{enumerate}
\end{thm}

Finally we have a generalization of the Green-Julg theorem, when $S_\red$ is unital and the
coaction $\delta_A$ is trivial. We will need the following notations to state it precisely:
\begin{itemize}  \setlength{\itemsep}{0ex}
\item[--] Let $(B\rtimes_\red\hat S)_1$ be the \Cst algebra $B\rtimes_\red\hat S$ endowed with
  the trivial coaction of $S_\red$. Taking into account the fact that any $B$-module $E$ can be
  fitted with the trivial coaction of $S_\red$, we get a morphism $\psi : KK(A,
  B\rtimes_\red\hat S) \to KK_{S_\red}(A, (B\rtimes_\red\hat S)_1)$.
\item[--] The inclusion $B\rtimes_\red\hat S \subset K(E)^{S_\red}$, where $E = B\tens H$ is
  endowed with the coaction $\delta'_E$, defines an element $\beta$ of
  $KK_{S_\red}({(B\rtimes_\red\hat S)_1}, B)$.
\item[--] In the compact case, the trivial representation factorizes through $\hat S_\red$ and
  its central support $p_0$ is an element of $\hat S_\red$. We call $\phi$ the homomorphism
  from $A$ to $A\rtimes_\red\hat S = A\tens \hat S_\red$ given by $\phi(a) = a\tens p_0$.
\end{itemize}
\begin{thm} \label{thm:green_julg}
  Let $S_\red$ be the reduced Hopf \Cst algebra associated to a compact quantum group ($S_\red$
  is unital). Assume that $S_\red$ coacts trivially on $A$, and that $B$ is a $S_\red$-algebra.
  Then the morphism
  \begin{displaymath}
    \Phi_1 : KK_{S_\red}(A,B) \verylongarrow{j_{\red}}
      KK(A\rtimes_\red\hat S, B\rtimes_\red\hat S) \verylongarrow{\phi^*}
      KK(A,B\rtimes_\red\hat S)
  \end{displaymath}
  is an isomorphism, and its inverse is given by
  \begin{displaymath}
    \Phi_2 : KK(A,B\rtimes_\red\hat S) \verylongarrow{\psi}
      KK_{S_\red}(A,(B\rtimes_\red\hat S)_1) \verylongarrow{\cdot\,\tens\beta}
      KK_{S_\red}(A,B) \text{.}
  \end{displaymath}
\end{thm}
As in the classical case, there is a ``dual statement'' of the Green-Julg theorem, which is
much easier to prove: if $S_\red$ is unital and if $\hat S_\red$ coacts trivially on $B$, there
is a canonical isomorphism $\Psi : KK_{\hat S_\red}(A,B) \to KK(A\rtimes S,B)$.

\bigskip

Let us close this section with the notion of $K$-amenability, which was introduced in
\cite{Cuntz:kmoy} for discrete groups and in \cite{JulgValette:kmoy} for locally compact
groups. We say that a locally compact quantum group $(\hat S_\red, \hat\delta_\red)$ is
$K$-amenable if the unit element of $KK_{\hat S_\red}(\CC,\CC)$ can be represented by a triple
$(E,\pi,F)$ such that the representation of the quantum group on $E$ is weakly contained in its
regular representation. Using the tools introduced in this section, most notably the descent
morphisms and the homomorphism $\delta'_A$, one can give other characterizations of
$K$-amenability, at least in the discrete case:
\begin{thm} \label{thm:K_amenability}
  Let $S$ and $S_\red$ be the full and reduced dual Hopf \Cst algebras of a locally compact
  quantum group $(\hat S_\red,\hat\delta_\red)$. Let $\varepsilon : S\to \CC$ be the trivial
  representation of $S$.  Then we have \ref{enum:kmoy_un} $\Rightarrow$ \ref{enum:kmoy_croise}
  $\Rightarrow$ \ref{enum:kmoy_lambda} $\Rightarrow$ \ref{enum:kmoy_alpha} and, if $S_\red$ is
  unital, \makebox{\ref{enum:kmoy_alpha} $\Rightarrow$ \ref{enum:kmoy_un}}:
  \begin{enumerate} \setlength{\itemsep}{0ex}
    \renewcommand\theenumi{{\it\roman{enumi}}} \renewcommand\labelenumi{{\it\roman{enumi}}.}
  \item \label{enum:kmoy_un} $(\hat S_\red, \hat\delta_\red)$ is $K$-amenable.
  \item \label{enum:kmoy_croise} For every \Cst algebra $A$ endowed with a coaction of $\hat
    S_\red$, \newline $[\lambda_A] \in$  ${KK(A\rtimes S, A\rtimes_\red S)}$ is invertible.
  \item \label{enum:kmoy_lambda} $[\lambda] \in KK(S, S_\red)$ is invertible.
  \item \label{enum:kmoy_alpha} There exists $\alpha \in KK(S_\red,\CC)$ such that
    $\lambda^*(\alpha) = [\varepsilon] \in KK(S,\CC)$.
  \end{enumerate}
\end{thm}

\section{Complements on discrete quantum groups}

We start with some useful results on subgroups and amalgamated free products of discrete
quantum groups. In particular, Proposition~\ref{prp:haar_libre_amalg} generalizes to the
amalgamated case the identity $h = h_1\free h_2$ for the Haar state of a free product
\cite{Wang:freeprod}. Together with Proposition~\ref{prp:sous_groupe_esp}, it shows that the
reduced Hopf \Cst algebra of an amalgamated free product $S_1\free_{T} S_2$ can be identified
with the reduced amalgamated free product of the reduced Hopf \Cst algebras.

\bigskip

In this section discrete quantum groups will be given by (one of) their \Cst algebra(s), i.e. by
a Woronowicz \Cst algebra $(S,\delta)$. Let us recall that a Woronowicz \Cst algebra is a
unital Hopf \Cst algebra such that the subspaces $(1\tens S)\delta(S)$ and $(S\tens
1)\delta(S)$ are dense in $S\tens S$. These conditions ensure the existence of a unique Haar
state on $S$, i.e. a state $h\in S'$ such that $(h\tens\id)\rond\delta = (\id\tens h)\rond\delta
= 1_Sh$. The GNS representation of $h$ is called the regular representation, and its image
$S_\red$ is the reduced Woronowicz \Cst algebra of the discrete quantum group. We say that $S$
is reduced if its regular representation is faithful, and in this case the Haar state $h$ is
also faithful (cf the proof of Proposition~\ref{prp:sous_groupe_esp}).

The general theory of Tannaka-Krein duality \cite{Woro:dual} shows that the discrete quantum
group associated to $(S,\delta)$ is fully characterized by the category $\mathcal{C}$ of the
finite-dimensional unitary corepresentations of $(S,\delta)$. We will use the notation
$\mathrm{Irr}~\mathcal{C}$ for a complete system of irreducible corepresentations in
$\mathcal{C}$. As a direct consequence of the general theory, we can see the subgroups of a
discrete quantum group in both pictures:

\begin{lem}\label{prp:sous_groupe_def}
  Let $(S,\delta)$ be a Woronowicz \Cst algebra and $\mathcal{C}$ the category of its
  finite-dimensional unitary corepresentations. There is then a natural bijection between:
  \begin{itemize}
  \item[--] Woronowicz sub-\Cst algebras $T\subset S$ (with $1_T = 1_S$),
  \item[--] full subcategories $\mathcal{D} \subset \mathcal{C}$ such that $1_{\mathcal{C}} \in
    \mathcal{D}$, $\mathcal{D}\tens \mathcal{D}\subset \mathcal{D}$ and $\bar{\mathcal{D}} =
    \mathcal{D}$.
  \end{itemize}
  Given such a subcategory $\mathcal{D}\subset \mathcal{C}$, the corresponding Woronowicz
  sub-\Cst algebra $T$ is the closed subspace of $S$ generated by the coefficients of the
  corepresentations $r\in \mathcal{D}$.
\end{lem}

\begin{prp}\label{prp:sous_groupe_esp}
  Let $(S,\delta)$ be a reduced Woronowicz \Cst algebra, and $T$ a Woronowicz sub-\Cst
  algebra of $S$. Then:
  \begin{enumerate}
  \item \label{item:prp_sous_groupe_haar} the Haar state $h_T$ is the restriction of $h_{S}$,
    and $T$ is reduced,
  \item \label{item:prp_sous_groupe_esp} there exists a unique conditional expectation $P : S
    \twoheadrightarrow T$ such that $h_S = h_T \rond P$,
  \item \label{item:prp_sous_groupe_inv} this conditional expectation is also characterized by
    the invariance property $(\id\tens P) \rond\delta_S = (P\tens\id) \rond\delta_S = \delta_T
    \rond P$.
  \end{enumerate}
\end{prp}

\begin{dem}  ~ \\ \indent
  \ref{item:prp_sous_groupe_haar}. The fact that $h_T$ is the restriction of $h_S$ to $T$
  follows immediately from the uniqueness of $h_T$. Let $\lambda_S$, $\lambda_T$ be the
  respective regular representations of $S$ and $T$. Because Haar states of quantum groups are
  KMS, we have
  \begin{eqnarray*}
    \mathrm{Ker}~\lambda_T &=& \{y\in T ~|~ h_T(y^*y)=0\} \\
    &=& \{y\in T ~|~ h_S(y^*y)=0\} ~=~ \mathrm{Ker}~\lambda_S \cap T ~=~ \{0\} \text{.}
  \end{eqnarray*}
  
  \ref{item:prp_sous_groupe_esp}.  Let $\mathcal{D} \subset \mathcal{C}$ be the categories
  of finite-dimensional unitary corepresentations of $T$ and $S$ respectively. In this
  paragraph we will repeatedly use the structure theorem of Woronowicz \cite{Woro:cmp}
  giving the Hopf \Cst algebra structure of $S$ in terms of the category $\mathcal{C}$.
  Let $\mathcal{S}^\circ$ (resp.  $\mathcal{T}$) be the linear span in $S$ of the
  coefficients of the corepresentations $r\in \mathrm{Irr}~\mathcal{C} \setminus
  \mathrm{Irr}~\mathcal{D}$ (resp.  $r\in\mathrm{Irr}~ \mathcal{D}$) : the direct sum
  $\mathcal{S} = \mathcal{S}^\circ \oplus \mathcal{T}$ is dense in $S$.  As a consequence
  of the stability properties of $\mathcal{D}$, we have the inclusions $\mathcal{S}^\circ
  \mathcal{T} \subset \mathcal{S}^\circ$, $\mathcal{T}\mathcal{S}^\circ \subset
  \mathcal{S}^\circ$ and $\mathcal{S}^{\circ *} \subset \mathcal{S}^\circ$. For instance
  if $s \in \mathrm{Irr}~\mathcal{C} \setminus \mathrm{Irr}~\mathcal{D}$, $t\in
  \mathrm{Irr}~\mathcal{D}$ and $r\subset s\tens t$ then $s\subset r\tens\bar t$, so that
  $r$ cannot be in $\mathcal{D}$. When passing to coefficients one obtains the first
  inclusion. In particular, these inclusions imply that $\Lambda_h(\mathcal{S}^\circ)$ and
  $\Lambda_h(\mathcal{T})$ are orthogonal subspaces in the GNS construction of $h_S$,
  because $h_S$ maps the coefficients of all irreducible corepresentations to $0$, except
  the ones of the trivial corepresentation.
  
  Now we identify $S$ and $T$ with their images in the regular representation of $S$. Let $p$
  be the orthogonal projection onto the closure of $\Lambda_h(\mathcal{T})$. The results of the
  previous paragraph (and the fact that $T$ is reduced) show that $Tp = pTp \simeq T$ and
  $p\mathcal{S}^\circ p = \{0\}$. Hence the compression by $p$ defines a conditional expectation
  $P : S \to pTp \simeq T$. The identity $h_S = h_T \rond P$ follows from the fact that the
  Haar states are the vectorial states associated to $\Lambda_h(1)$. Conversely, if $P' : S \to
  T$ is a conditional expectation such that $h_S = h_T \rond P'$, for every $x\in S$ we have
  $P'(x) = P(x)$ because $T$ is reduced and
  \begin{displaymath}
    \forall~ y\in T~~ h_T(P'(x)y) = h_T(P'(xy)) = h_S(xy) = h_T(P(xy)) = h_T(P(x)y) \text{.}
  \end{displaymath}
  
  \ref{item:prp_sous_groupe_inv}. Finally, the invariance property of $P$ can be easily
  verified on $\mathcal{T}$ and $\mathcal{S}^\circ$ separately: we have $\delta(\mathcal{T})
  \subset \mathcal{T}\tens\mathcal{T}$ and $\delta(\mathcal{S}^\circ) \subset
  \mathcal{S}^\circ\tens\mathcal{S}^\circ$, but $P$ coincides with $\id$ on $T$, and with $0$
  on $\mathcal{S}^\circ$. Conversely, if $P' : S \to T$ is a conditional expectation verifying
  the invariance property, we have
  \begin{displaymath}
    1_T h_S = (\id\tens h_S) (P'\tens\id) \delta_{S} = (\id\tens h_T) \delta_T P' = 1_Th_T \rond P'
  \end{displaymath}
  so that $h_S = h_T\rond P'$ and $P' = P$.
\end{dem}

\begin{lem} \label{lem:cat_quotient}
  Let $(S,\delta)$ be a reduced Woronowicz \Cst algebra, and $T$ a \linebreak Woronowicz
  sub-\Cst algebra of $S$. Let $E$ be the Hilbert $T$-module of the GNS representation of $S$
  associated to the conditional expectation $P$ introduced in
  Proposition~\ref{prp:sous_groupe_esp}.
  \begin{enumerate}
  \item \label{enum:lem_cat_quotient_def} The relation on $\mathrm{Irr}~\mathcal{C}$ defined by
    $r \sim r' \Leftrightarrow (\exists~ t\in\mathcal{D}~~ t\subset \bar r\tens r')$ is an
    equivalence relation, let $(\mathrm{Irr}~\mathcal{C}) / \mathcal{D}$ be the corresponding
    quotient set.
  \item \label{enum:lem_cat_quotient_compl} Let $\mathcal{S}^r$ be the subspace of
    $\mathcal{S}$ generated by the coefficients of a corepresentation $r\in \mathcal{C}$. For
    $\alpha \in (\mathrm{Irr}~\mathcal{C}) / \mathcal{D}$ we put $\mathcal{S}^\alpha =
    \sum_{r\in\alpha} \mathcal{S}^r$, let $S^\circ$ be the closure of $\mathcal{S}^\circ =
    \sum_{\alpha\neq 1} \mathcal{S}^\alpha$. Then we have $\mathrm{Ker}~P = S^\circ$.
  \item \label{enum:lem_cat_quotient_decomp} We denote by $E^\alpha$ and $E^r$ the respective
    closures of $\Lambda_P(\mathcal{S}^\alpha)$ and $\Lambda_P(\mathcal{S}^r)$ in $E$. For $r
    \in \alpha$, the closed submodule of $E$ generated by $E^r$ is $E^\alpha$. Moreover $E$ is
    the orthogonal direct sum of the $E^\alpha$.
  \end{enumerate}
\end{lem}

\begin{dem} ~ \\ \indent
  \ref{enum:lem_cat_quotient_def}. The first point essentially results from the following
  equivalence, that we already used in the previous proof: $t\subset \bar r\tens r'
  \Leftrightarrow r'\subset r\tens t$. We therefore have $r\sim r' \Leftrightarrow r'\subset
  r\tens\mathcal{D}$, and the transitivity of $\sim$ follows from the hypothesis
  $\mathcal{D}\tens \mathcal{D} \subset \mathcal{D}$. The other conditions are obvious on the
  definition (and result from the hypothesis $\bar{\mathcal{D}} = \mathcal{D}$,
  $1_{\mathcal{C}} \in \mathcal{D}$).
  
  \ref{enum:lem_cat_quotient_compl}. It is clear that the class of $1_{\mathcal{C}}$ in
  $(\mathrm{Irr}~\mathcal{C}) / \mathcal{D}$ is $\mathrm{Irr}~\mathcal{D}$. In particular the
  subspace $\mathcal{S}^\circ$ coincides with the one we have used in the proof of
  Proposition~\ref{prp:sous_groupe_esp}, and we have already seen that it is included in
  $\mathrm{Ker}~ P$. Now if we write $x\in \mathrm{Ker}~P$ as a limit of sums $y_n+z_n$ with
  $y_n \in \mathcal{T}$ and $z_n \in \mathcal{S}^\circ$, we have $P(y_n + z_n) = y_n
  \rightarrow 0$, so that $x$ is in the closure of $\mathcal{S}^\circ$.
  
  \ref{enum:lem_cat_quotient_decomp}. The linear span of $\mathcal{S}^r\mathcal{T}$ coincides
  with the subspace generated by the coefficients of the corepresentations $r\tens t$,
  $t\in\mathcal{D}$, which is exactly $\mathcal{S}^\alpha$: recall that $r'\sim r
  \Leftrightarrow r'\subset r\tens\mathcal{D}$. So the closed submodule of $E$ generated by
  $E^r$ is $E^\alpha$. Now let $\alpha$ and $\beta$ be two distinct elements of
  $(\mathrm{Irr}~\mathcal{C}) / \mathcal{D}$, and take $x \in \mathcal{S}^\alpha$, $y\in
  \mathcal{S}^\beta$. We can assume that $x$ and $y$ are respective coefficients of two
  corepresentations $r\in\alpha$, $s\in\beta$. The subobjects of $\bar r\tens s$ are not in
  $\mathcal{D}$ because $\alpha$ and $\beta$ are distinct, so that $x^*y \in \mathcal{S}^\circ$
  and
  \begin{displaymath}
    \left<\Lambda_P(x),\Lambda_P(y)\right> = P(x^*y) = 0 \text{.}
  \end{displaymath}
  This establishes the fact that $E^\alpha$ and $E^\beta$ are orthogonal, and it is clear that
  the sum of the subspaces $E^\alpha$ is dense in $E$.
\end{dem}

We will use the following notations about reduced amalgamated free products
\cite{Avitzour:free,Voiculescu:freesymm,VoicuDykemaNica:freerandom}. Let $S_1$, $S_2$ be unital
\Cst algebras, and $T$ a common sub-\Cst algebra of $S_1$ and $S_2$ containing their unit
element. Let $P_1$, $P_2$ be two conditional expectations of $S_1$ and $S_2$ onto $T$, and
$(E_1,\eta_1)$, $(E_2,\eta_2)$ the corresponding GNS constructions. Like in
Lemma~\ref{lem:cat_quotient}, we denote by $S_i^\circ$ the kernel of $P_i$, and by $E_i^\circ$
the orthogonal of $\eta_i$ in $E_i$ --- i.e. the closure of $\Lambda_{P_i}(S_i^\circ)$.
Following \cite{Voiculescu:freesymm}, we put \renewcommand{\theequation}{$*$}
\begin{eqnarray} \label{eq:free_prod_mod}
  & E_1\free_T E_2 = \eta T \oplus \textstyle\bigoplus_{n=1}^\infty \bigoplus_{(i_j)\in I_n}
  (E_{i_1}^\circ\tens_T \cdots \tens_T E_{i_n}^\circ) \text{,~~~where} \\ \nonumber & I_n = \{
  (i_1,\ldots,i_n) \in \{1,2\}^n ~~|~~ \forall~k~~i_k\neq i_{k+1} \} \text{.}
\end{eqnarray}
The Hilbert $T$-module $E_1\free_T E_2$ carries then a natural representation of the (full)
free product with amalgamation $S_1\free_T S_2$, whose image is called the reduced free product
with amalgamation $(S_1,P_1) \free_T (S_2,P_2)$. The vacuum vector $\eta$ defines a conditional
expectation $P_1\free_T P_2$ from $S_1\free_T S_2$ onto $T$. Finally, we will denote by
$E(r,i)$ the closed submodule of $E$ associated to the $n$-tuples $(i_j)\in I_n$
in~(\ref{eq:free_prod_mod}) such that $i_n \neq i$.

In \cite{Wang:freeprod}, it is shown that the free product of two \mbox{Woronowicz} \Cst
algebras with amalgamation over a Woronowicz sub-\Cst algebra carries again a natural structure
of Woronowicz \Cst algebra. Moreover the Haar state and the category of finite-dimensional
unitary corepresentations of that free product are investigated in the non-amalgamated case. It
is not hard to generalize the result on the Haar state to the amalgamated case:

\begin{prp} \label{prp:haar_libre_amalg}
  Let $S'_1$, $S'_2$ be two Woronowicz \Cst algebras, and $T'$ be a common Woronowicz sub-\Cst
  algebra of $S'_1$ and $S'_2$. Let $S_1$, $S_2$ and $T$ be their respective reduced Woronowicz
  \Cst algebras. We consider the free product with amalgamation $S = S_1\free_T S_2$ and denote
  by $\lambda : S \to S_\red$ the reduction homomorphism of $(S,\delta)$.
  \begin{enumerate}
  \item \label{enum:prp_haar_amalg} If $P$, $P_1$ and $P_2$ are the canonical conditional
    expectations of $S_\red$, $S_1$ and $S_2$ onto $T$, we have $P_1\free_T P_2 =
    P\rond\lambda$.
  \item \label{enum:prp_red_amalg} The reduced Woronowicz \Cst algebra of $S'_1\free_{T'} S'_2$
    is isomorphic to the reduced amalgamated free product of $(S_1, P_1)$ and $(S_2, P_2)$.
  \end{enumerate}
\end{prp}

\begin{dem} ~ \\ \indent
  \ref{enum:prp_haar_amalg}. We put $P' = P_1\free_T P_2$, and we prove in the next paragraph
  that $h_T\rond P'$ is the Haar state of $S$. This will prove that $P'$ factors through
  $S_\red$ and, by Proposition~\ref{prp:sous_groupe_esp}, that its factorization equals $P$.
  Thanks to the uniqueness of the Haar state, it is enough to establish the invariance of the
  state $h_T\rond P'$.
  
  The \Cst algebra $S$ is a quotient of the free product $S_1\free_\CC S_2$, so that the
  elements $s_1 \cdots s_n$ with $s_k \in \mathcal{S}_{i_k}$ and $(i_j) \in I_n$ generate the
  normed vector space $S$. Moreover it is easily seen by induction, using the decompositions
  $\mathcal{S}_i = \mathcal{T} \oplus \mathcal{S}_i^\circ$, that it is enough to consider the
  elements of $\mathcal{T}$ and the elements $s_1 \cdots s_n$ with $s_k \in
  \mathcal{S}_{i_k}^\circ$ and $(i_j) \in I_n$. On the subspace $\mathcal{T}$, the expectation
  $P'$ coincides with the identity, so that the invariance property is evident.
  
  Therefore we consider an element $s = s_1 \cdots s_n$ with $s_k \in \mathcal{S}_{i_k}^\circ$
  and $(i_j)\in I_n$. For such an element we have $P'(s) = 0$ because $S_1$, $S_2$ are free in
  $(S,P')$ and $\mathcal{S}_i^\circ \subset \mathrm{Ker}~ P_i$. On the other hand, each
  $\delta_{i_k}(s_k)$ can be written as a sum $\sum_l s_{k,l}\tens s'_{k,l}$ with $s_{k,l}$,
  $s'_{k,l} \in \mathcal{S}_{i_k}^\circ$, and one obtains for $\delta(s) = \delta_{i_1}(s_1)
  \cdots \delta_{i_n}(s_n)$ an expression of the form $\sum_l s_l\tens s'_l$, where the $s_l$
  and $s'_l$ are all in $\mathrm{Ker}~ P'$ --- using again the freeness of $S_1$ and $S_2$.
  Therefore we also have $(P'\tens\id)\rond\delta (s) = 0 $ and $(\id\tens P')\rond\delta (s) =
  0$.
  
  \ref{enum:prp_red_amalg}. The Hopf \Cst algebras $S_1\free_T S_2$ and $S'_1\free_{T'} S'_2$
  having the same dense Hopf sub-$*$-algebras of coefficients, they have the same reduced
  Hopf \Cst algebras. Now the first point shows that the Haar state of $S_1\free_T S_2$ is
  $h_T\rond (P_1\free_T P_2)$. But $T$ is reduced, so that the reduced Hopf \Cst algebra of
  $S_1\free_T S_2$ is exactly $(S_1,P_1) \free_{T} (S_2,P_2)$.
\end{dem}

\begin{rque}{Remark}
  In the non-amalgamated case, it is quite easy to compute the category of finite-dimensional
  unitary corepresentations of $S$: the corepresentations $v_1\tens\cdots\tens v_n$ , where
  $v_k\in \mathrm{Irr}~\mathcal{C}_{i_k} \setminus \{1\}$ and $(i_k) \in I_n$, form a complete
  system of irreducible corepresentations of $S$, and the fusion rules can be naturally deduced
  from the ones of $\mathcal{C}_1$ and $\mathcal{C}_2$ \cite{Wang:freeprod}.  This is very
  similar to the structure theorem for the free product of two discrete groups.  
  
  There is no such ``simple'' generalization of the classical theory in the amalgamated case.
  More precisely, if $v_1$ and $v_2$ are respective irreducible corepresentations of $S_1$ and
  $S_2$, but not of $T$, the corepresentation $v_1\tens v_2$ does not need to be irreducible for
  $S$: see the example hereafter. However, looking at the amalgamated free product modulo $T$,
  one recovers a part of the classical structure: this is the meaning of the decomposition of
  $E_1\free_T E_2$ into orthogonal $T$-modules given by Lemma~\ref{lem:cat_quotient}
  and~(\ref{eq:free_prod_mod}).
\end{rque}

\begin{rque}{Example}
  Consider two copies $S_1$, $S_2$ of $A_o(Q)$, the Woronowicz \Cst algebra of an orthogonal
  free quantum group \cite{Wang:freeprod,DaeleWang:univ,Banic:O(n)_cras}. Let us recall that
  $\mathrm{Irr}~\mathcal{C}_1$ and $\mathrm{Irr}~\mathcal{C}_2$ can be indexed by $\NN$, with
  the same fusion rules as the ones of $SU(2)$: we will denote by $v_{i,k}$ the corresponding
  irreducible corepresentations, with $i\in \{1,2\}$ and $k\in\NN$. The discrete quantum group
  corresponding to $A_o(Q)$ has a unique non-trivial subgroup, associated to the subcategory
  generated by $\mathrm{Irr}~ \mathcal{D} = \{v_{2k} ~|~ k\in\NN\}$. Its Woronowicz \Cst
  algebra $T$ is the sub-\Cst algebra of even elements for the natural $\ZZ/2\ZZ$-grading of
  $A_o(Q)$ --- when $Q = I_n$, $T$ is also the Woronowicz \Cst algebra $A_{\mathrm{aut}}(M_n)$
  of \cite{Banic:symmetries}.
  
  The tensor product of the fundamental corepresentations of $S_1$ and $S_2$, $v_{1,1}\tens
  v_{2,1}$, defines a corepresentation of the amalgamated free product $S = S_1\free_T S_2$,
  and it is not hard to see that it admits a strict sub-corepresentation $a$ of dimension $1$.
  More generally the sub-corepresentations of $v_{i_1,k_1} \tens \cdots \tens v_{i_n,k_n}$,
  with $(i_l) \in \{1,2\}^n$ and $(k_l) \in \NN^n$, have the same dimensions and multiplicities as
  the ones of $v_{1,k_1}\tens\cdots\tens v_{1,k_n}$. Moreover one can prove, using the freeness
  of $S_1$ and $S_2$ in $S$, that $a$ generates a copy of $C^*(\ZZ)$ in $S$. Finally the
  corepresentations $a^k \tens v_{1,l}$, with $k\in\ZZ$ and $l\in\NN$, form a complete system
  of irreducible corepresentations of $S$. The fusions rules are induced by the ones of
  $C^*(\ZZ)$ and $A_o(Q)$, plus the relation $v_{1,1}\tens a \simeq v_{2,1} \simeq a^{-1}\tens
  v_{1,1}$.
\end{rque}

\section{The Bass-Serre tree and $K$-amenability}

The aim of this section is to prove the $K$-amenability of amalgamated free products of
amenable discrete quantum groups. This extends to the quantum case a result of Julg and Valette
\cite[cor.~4.1]{JulgValette:kmoy}, which was obtained as a corollary of their more general
result about locally compact groups acting on trees. Let us note however that the homotopy we
construct in Theorem~\ref{thm:amalg_kmoy} is closer to the one used by Cuntz for the case of
free products $G\free\ZZ$ \cite{Cuntz:kmoy}. Indeed the homotopy of Julg and Valette relies on
the construction of negative type functions on trees, which seems more complicated to extend to
the quantum case.

\bigskip

Let $S = S_1\free_T S_2$ be an amalgamated free product of reduced Woronowicz \Cst algebras.
Let $P$, $R_i$ be the canonical conditional expectations of $S_\red$ onto $T$, $S_i$
respectively, and $E$, $F_i$ the associated GNS constructions. For the definition of the
Julg-Valette operators we will need a slight reinforcement of the well-known isomorphisms
$E(r,i)\tens_T E_i \simeq E$:

\begin{lem} \label{lem:right_coset}
  The map $\Psi_i : \Lambda_P(x)\tens y \mapsto \Lambda_{R_i}(xy)$ defines an isomorphism
  between $E(r,i)\tens_T S_i$ and $F_i$.
\end{lem}

\begin{dem}
  We take $i=1$ for the proof. In the definition of $\Psi_1$ one can take $x\in\mathcal{T}$ or
  $x = y_1 \cdots y_n$, $y_k \in \mathcal{S}_{i_k}^\circ$, $(i_k) \in I_n$ and $i_n=2$. In
  particular the surjectivity of $\Psi_1$ will follow immediately if we prove that it is an
  isometry. So we have to show that $y^*P(x^*x)y = R_1(y^*x^*xy)$ for $y\in S_1$ and $x$ like
  above. By $S_1$-linearity of $R_1$ we can assume that $y=1$, and then the equality $P(x^*x) =
  R_1(x^*x)$ amounts to the fact that $R_1(x^*x) \in \mathcal{T}$, which is clear when $x\in
  \mathcal{T}$.
  
  So we take $x = y_1 \cdots y_n$ with $y_k \in \mathcal{S}^\circ_{i_k}$, $(i_k)\in I_n$,
  $i_n=2$, and we prove that $R_1(x^*x)$ is orthogonal to $\mathcal{S}_1^\circ$ in the GNS
  construction of $P$. For any $x_1\in \mathcal{S}^\circ_1$, we have
  \begin{displaymath}
    P(x_1 R_1(x^*x)) = P(R_1(x_1x^*x)) = P(x_1x^*x)
  \end{displaymath}
  by $S_1$-linearity of $R_1$. To conclude we proceed by induction on $n$ in the expression $x
  = y_1 \cdots y_n$ --- we put $n=0$ when $x\in\mathcal{T}$. By
  Proposition~\ref{prp:haar_libre_amalg}, $S_1$ and $S_2$ are free in $S_\red$ with respect to
  $P$. Therefore
  \begin{eqnarray*}
    P(x_1x^*x) &=& P(x_1y_n^*\cdots y_2^*y_1^*y_1y_2 \cdots y_n) \\
    &=& P(x_1y_n^*\cdots y_2^* P(y_1^*y_1)y_2 \cdots y_n) \\
    &=& P(x_1y_n^*\cdots y^{\prime *}_2 y'_2\cdots y_n)
  \end{eqnarray*}
  with $y'_2 = P(y_1^*y_1)^{1/2}y_2$, and the last expression equals $0$ by induction.
\end{dem}

We now assume that $S_1$ and $S_2$ are amenable Woronowicz \Cst algebras. This means that they
are reduced and admit co-units (also called trivial representations in our context), i.e.
continuous characters $\varepsilon : S_i\to \CC$ such that $(\id\tens\varepsilon)\rond \delta =
(\varepsilon\tens\id) \rond\delta = \id$. In the following definition, we use the isomorphism
$F_i\tens_\varepsilon\CC \simeq E(r,i)\tens_\varepsilon\CC$ induced by $\Psi_i$.

\begin{df} \label{df:graphe_amalg}
  Let $S = S_1\free_T S_2$ be an amalgamated free product of two amenable Woronowicz \Cst
  algebras. Let $(E,\eta)$, $(F_i,\eta_i)$ be the GNS constructions associated to the canonical
  conditional expectations of $S$ onto $T$, $S_i$ respectively. Let $\varepsilon$ denote the
  co-unit of $S_1$ and $S_2$.
  \begin{enumerate}
  \item The Hilbert space $K_g = E\tens_\varepsilon \CC$ (resp. $H = F_1\tens_\varepsilon\CC
    \oplus F_2\tens_\varepsilon\CC$) is called the $\ell^2$-space of geometric edges (resp. of
    vertices) of the quantum graph associated to $(S_1, S_2, T)$. They carry the GNS
    representations of $S$.
  \item For that quantum graph, we define two Julg-Valette operators $\Phi_i$, \linebreak
    ${i\in\{1,2\}}$, given by the following formulas (where $\{i,j\} = \{1,2\}$) :
    \begin{displaymath}
      \Phi_i : \left\{ \begin{array}{l}
          H \to K_g \text{,} \\
          \eta_j\tens_\varepsilon 1 \mapsto \eta\tens_\varepsilon 1 \text{,~~~}
          \eta_i\tens_\varepsilon 1 \mapsto 0 \text{,} \\
          F_k^\circ\tens_\varepsilon\CC \stackrel {\sim} {\longrightarrow}
          E(r,k)^\circ\tens_\varepsilon\CC \text{~~~ for $k=1,2$.}
          \end{array} \right.
    \end{displaymath}
  \end{enumerate}
\end{df}

\begin{rque}{Example}
  When $S_1 = C^*_\red(G_1)$, $S_2 = C^*_\red(G_2)$ and $T = C^*_\red(H)$, the \Cst
  algebra $S$ is some \Cst completion of $\ell^1(G)$, with $G = G_1\free_H G_2$. Then $K_g
  \simeq \ell^2 (G/H)$ and $H \simeq \ell^2(G/G_1 \sqcup G/G_2)$ are the $\ell^2$-spaces
  of edges and vertices of the Bass-Serre tree associated to $(G_1, G_2, H)$
  \cite{Serres:arbres}. Note that the isomorphism $F_i\tens_\varepsilon\CC \simeq
  E(r,i)\tens_\varepsilon\CC$ in this case is just the $\ell^2$-formulation of the
  canonical bijection $G(r,i)/H \stackrel {\sim} {\longrightarrow} G/G_i$, where $G(r,i)
  \subset G$ is the subset of alternate products $g_{i_1}\cdots g_{i_k}$ that ``do not end
  in $G_i$'' (including $H$).
  
  In the Bass-Serre tree, there are two evident ways of choosing one endpoint for each
  edge. They are given by the following maps from $G/H$ to $G/G_1 \sqcup G/G_2$, that
  determine the graph structure:
  \begin{displaymath}
    xH \mapsto xG_1 \text{~~~and~~~} xH \mapsto xG_2 \text{.}
  \end{displaymath}
  However, the corresponding ``target operators'' on the $\ell^2$-spaces do not need to be
  bounded. If we choose $G_i \in G/G_i$ as the origin of the tree, there is another natural,
  more geometrical way of selecting one endpoint for each edge: namely the furthest one from
  the origin. The adjoint of the corresponding target operator is then precisely the
  Julg-Valette operator $\Phi_i$ of Definition~\ref{df:graphe_amalg}.
\end{rque}

\begin{thm} \label{thm:amalg_kmoy}
  Let $S = S_1\free_T S_2$ be an amalgamated free product of two amenable Woronowicz \Cst
  algebras. Take $i\in \{1,2\}$.
  \begin{enumerate}
  \item \label{item:thm_amalg_reg} The GNS representations of $S$ on the Hilbert spaces $H$ and
    $K_g$ introduced in Definition~\ref{df:graphe_amalg} are weakly contained in the regular
    representation.
  \item \label{item:thm_amalg_kk} The Julg-Valette operator $\Phi_i : H\to K_g$ is a Fredholm
    operator and commutes to the representations of $S$ modulo compact operators. Let $\alpha$
    be the associated element of $KK(S_\red,\CC)$.
  \item \label{item:thm_amalg_kmoy} The element $\lambda^*(\alpha)\in KK(S,\CC)$ is equal to
    the element $[\varepsilon]\in KK(S,\CC)$ associated to the co-unit of $S$. Hence the
    discrete quantum group associated to $S$ is $K$-amenable.
  \end{enumerate}
\end{thm}

\begin{dem} ~ \\ \indent
  \ref{item:thm_amalg_reg}. The Woronowicz \Cst algebras $S_1$ and $S_2$ being amenable, so is
  $T$, and this means that its trivial representation is weakly contained in the regular one.
  As a consequence, the GNS representation of $\varepsilon\rond P$ is weakly contained in the
  GNS representation of $h \rond P$. But these representations are respectively the GNS
  representation on $E\tens_\varepsilon\CC$ and the regular representation, by
  Proposition~\ref{prp:sous_groupe_esp}. The same argument works for $F_1$ and $F_2$.
  
  \ref{item:thm_amalg_kk}. It is clear on the definition, and thanks to
  Lemma~\ref{lem:right_coset}, that the Julg-Valette operators $\Phi_i$ are surjective and have
  a kernel of dimension one. Hence they are Fredholm operators. We will now prove that $\Phi_1$
  commutes to the action of an element $s \in S_1$. For symmetry reasons, this will imply that
  $\Phi_2$ commutes to $S_2$, and because $\Phi_1$ and $\Phi_2$ are equal modulo a compact
  operator, they will both commute to $S_1 \cup S_2$, hence to $S$, modulo compact operators.
  
  The main point is the following one: if $\zeta$ is a vector of $E(r,k)$ such that $s\zeta$ is
  again in $E(r,k)$, we have evidently
  \begin{displaymath}
    \Psi_k(s\zeta\tens_T x) = s\Psi_k(\zeta\tens_T x)
  \end{displaymath}
  for every $x\in S_k$. But it is easy to check that $s E(r,1)^\circ \subset E(r,1)^\circ$ and
  $s E(r,2) \subset E(r,2)$ for $s\in S_1$. As $\Phi_1$ coincides respectively on
  $E(r,1)^\circ$ and $E(r,2)$ with the inverses of $\Psi_1\tens_\varepsilon\id$ and
  $\Psi_2\tens_\varepsilon\id$, it commutes with the action of $s$ on $E(r,1)^\circ \oplus
  E(r,2) = (\eta_1\tens_\varepsilon 1)^\bot$. Finally, $S_1$ acts trivially on the line $\eta_1
  \tens_\varepsilon \CC$, that $\Phi_1$ maps to $0$.
  
  \ref{item:thm_amalg_kmoy}.  Let $\tilde\Phi_i : H \to K_g \oplus \CC$ be the operator which
  coincides with $\Phi_i$ on $(\eta_i\tens_\varepsilon 1)^\bot$, but maps
  $(\eta_i\tens_\varepsilon 1)$ to $1_\CC$. Let us denote by $\pi$ the sum of the GNS and the
  trivial representations of $S$ on $K_g \oplus \CC$, and by $\rho$ the GNS representation on
  $H$. It is clear that $\tilde \Phi_i$ is a compact perturbation of $\Phi_i + 0 : H\to K_g
  \oplus \CC$, so that the class of $x = (\pi, \rho, \tilde\Phi_i)$ in $KK(S, \CC)$ equals
  $\lambda^*(\alpha) - [\varepsilon]$. We will now prove that $x$ is homotopic to a degenerate
  triple.
  
  It is clear from the proof of point~\ref{item:thm_amalg_kk} that $\tilde\Phi_i$ intertwines
  the restrictions of $\pi$ and $\rho$ to $S_i$. Moreover one has $\tilde\Phi_1 = \tilde\Phi_2
  \rond u$, where $u \in L(H)$ equals the identity on the subspaces
  $F_k^\circ\tens_\varepsilon\CC$ and exchanges $\eta_1\tens_\varepsilon 1$ and
  $\eta_2\tens_\varepsilon 1$. In particular, $\tilde\Phi_1$ intertwines $\rho(s)$ and $\pi(s)$ if
  $s\in S_1$, but $u\rho(s)u$ and $\pi(s)$ if $s\in S_2$.
  
  Now $u$ is unitary and self-adjoint, so that the expression $u_t = \cos t + {\rm i} u \sin t$
  defines a family of unitaries $(u_t)$. Because $u$ commutes to $\rho(T)$, one can define a
  family of representations of $S = S_1\free_T S_2$ on $H$ by the formulas $\rho_t(s) = u_t^*
  \rho(s) u_t$ for $s\in S_2$ and $\rho_t(s) = \rho(s)$ for $s\in S_1$. Leaving $\pi$ and
  $\tilde\Phi_1$ unchanged, this gives a homotopy between $x$ and a degenerate triple: by
  construction, $\tilde\Phi_1$ intertwines exactly $\rho_1$ and $\pi$ and is unitary.
\end{dem}

\begin{rque}{Remark}
  In \cite{Germain:freereduced2} Germain defines a topological property called ``relative
  $K$-nuclearity'' for inclusions of unital \Cst algebras $T\subset S_1$, $T\subset S_2$. This
  property mostly serves as a sufficient condition for the full free product with amalgamation
  $S = S_1\free_T S_2$ to be ``dominated in $K$-theory'' by the reduced free product, which is
  equivalent to $K$-amenability in the case of \Cst algebras of discrete quantum groups.
  However, when $T$ is infinite-dimensional very little is known about relative $K$-nuclearity,
  even in our case of inclusions of \Cst algebras of amenable discrete quantum groups. In this
  paper we chose a different, more geometrical way to prove the $K$-amenability of our
  amalgamated free products.
\end{rque}

\noindent {\sc Acknowledgments}

The results of this paper are part of my Ph.D. Thesis at the University Paris~7. I would
like to thank my supervisor G.~Skandalis for his numerous suggestions that helped to
simplify the proofs and the exposition. I it also a pleasure to thank P.~Julg for fruitful
discussions on the trees of amalgamated free products.

\bibliographystyle{plain} 
\bibliography{free}

\end{document}